# Characterizations of the Extended Geometric, Harris, Negative Binomial and Gamma Distributions

**Sandhya E [1], S. Sherly [2], Jos M. K [3] and Raju N [4]**

(1) Department of Statistics, Prajyoti Niketan College, Pudukkad, Trichur-680301, India.
 (*e-mail: esandhya@hotmail.com*)
(2) Department of Statistics, Vimala College, Trichur-680 009, Kerala, India.
 (*e-mail: sebastian_sherly@yahoo.com*)
(3) Department of Statistics, St. Thomas College, Trichur-680 001, Kerala, India.
 (*e-mail: jos_kuriakku@yahoo.com*)
(4) Department of Statistics, University of Calicut, Calicut University PO-673635, India.
 (*e-mail: na_raju@yahoo.com*)

**Abstract:** Extended geometric distribution is defined and its mixture is characterized by the property of having completely monotone probability sequence. Also, convolution equations and probability generating functions are used to characterize extended geometric distributions. Further, some characterizations of Harris and negative binomial distributions based on probability generating functions are obtained. Relations between these distributions are derived and finally a gamma distribution is characterized in terms of Laplace transform.

## 1. Introduction

A number of papers have been devoted to various characterizations of geometric distribution. A geometric distribution may be extended to cover the case of a variable taking values $a, a + k, a + 2k, \ldots (k>0)$. See, Johnson *et al*. (1992, p.201). In this paper, we consider a geometric distribution with parameter $p$ defined on $\{a, a + k, a + 2k, \ldots\}$ where $a \geq 0$ and $k > 0$ are integers. We call the corresponding distribution as extended geometric distribution and denote it by $Geo_a(p, k)$. In the notation the suffix $a$ suggests that the support of the distribution starts from $a$, and $k$ implies that the atoms (probability carrying integers) of the distribution are $k$ integers apart.

Let the random variable (r.v) $Z$ have a geometric distribution on $\{0, 1, 2, \ldots\}$ with parameter $p$ ( $Z \sim Geo(p)$ ) and probability mass function (p.m.f)

$$P(Z = n) = q^n p, \quad n = 0,1,2,\ldots, \quad 0 < p < 1, \quad p + q = 1.$$

---

The second author is grateful to the University Grants Commission of India for the Award of Teacher Fellowship under Faculty Improvement Programme.

Then the r.v $Y = a + kZ$ has an extended geometric distribution with p.m.f

$$P(Y = a + nk) = q^n p, \quad n = 0,1,2,\ldots \tag{1}$$

where $a \geq 0$ and $k > 0$ are integers, $0 < p < 1$ and $p + q = 1$. i.e., $Y \sim Geo_a(p, k)$.

The moments of an extended geometric distribution are easily available from the following simple but an important characterization property of the distribution.

A r.v $Z \sim Geo(p)$ if and only if the r.v $Y \sim Geo_a(p, k)$ where $Y = a + kZ$. Then,

$E(Y) = a + k\, E(Z)$

$V(Y) = k^2\, V(Z)$

Also the probability generating function (p.g.f) of $Geo_a(p, k)$ is of the form

$$P(s) = \frac{ps^a}{1 - qs^k}$$

Satheesh and Sandhya (1997) have given a characterization for mixtures of geometric distributions using the property of having completely monotone probability sequence (CMPS). They also characterized geometric distribution in terms of convolution equations. Similar characterizations developed in the case of extended geometric distributions are presented in section 2 of this paper.

Harris (1948) introduced a p.g.f

$$P(s) = \frac{s}{\left(m - (m-1)s^k\right)^{1/k}}, \quad k > 0 \text{ integer}, \ m > 1. \tag{2}$$

He discussed this p.g.f while considering a simple discrete branching process where a particle either splits into $(k+1)$ identical particles or remains the same during a short time interval $\Delta t$. The probability distribution corresponding to the above p.g.f is called Harris distribution and its properties are studied by Sandhya *et al.*(2005, available at *http://arxiv.org/abs/math.ST/0506220* ).

A r.v $X$ follows a Harris distribution with parameters $m$ and $k$ if its p.m.f is given by

$$f(x) = \binom{(1/k)+n-1}{n}\left(\frac{1}{m}\right)^{1/k}\left(1 - \frac{1}{m}\right)^n, \quad x = 1 + nk, \ n = 0,1,2,\ldots, \tag{3}$$

where $k>0$ is an integer and $m>1$. We write $X \sim H_1(m,k,1/k)$. It may be noted that unlike the well known standard discrete distributions, Harris distribution is concentrated on $\{1, 1+k, 1+2k, \ldots\}$, where $k$ is fixed. Characterizations of Harris and negative binomial distributions based on the p.g.f s are given in section 3. Also some relations between extended geometric distribution and Harris distribution are discussed in this section.

Relations among the Harris, geometric, negative binomial and gamma distributions are given in Sandhya *et al.* (2005). It is shown that Harris distribution is a generalization of decapitated geometric distribution to which it reduces when $k = 1$. Also it is proved that a r.v $X \sim H_1(m,k,1/k)$ if and only if the r.v $W = (X-1)/k$ follows a negative binomial distribution $NB(1/m, 1/k)$ with parameters $1/m$ and $1/k$. For a negative binomial distribution with index parameter $1/k$, see Cadigan and Chen (2001). The r.v $W$ denotes the number of failures preceding a fractional success. There are situations where a success can occur only when $k$ (>1, an integer) fractional successes happen. For example, suppose a sales representative receives an incentive only when he sells $k$ items. Such an example is considered in Sherly *et al.*(2005, available at *http://arxiv.org/abs/math.ST/0510658*). Here selling of each item is a fractional success and obtaining an incentive is a success and success occurs only when the $k$ fractional successes happen *i.e.* when $k$ items are sold. Let $X_i$ denote the number of failures preceding the $i^{th}$ fractional success, $i = 1, 2, \ldots, k$ and $X_i \sim NB(p,1/k)$ where the probability of success $p$ is assumed to be a constant for each Bernoulli trial. Then their sum $Y = \sum_{i=1}^{k} X_i$, the total number of failures preceding a success is a geometric r.v. But a well known result connecting the geometric and negative binomial distribution is that the sum of $k$ independently and identically distributed (i.i.d) geometric r.v s with parameter $p$ is a $NB(p, k)$ r.v. See Johnson *et al.*(1992, p.203). A characterization of negative binomial distribution $NB(p,1/k)$ is given in section 3 of this paper.

The following are some other similar real life examples:

(a) In product control, an item is passed a quality test only when $k$ identical quality tests are passed.

(b) The strong room of a bank is broken only when $k$ identical locks are broken.

(c) A course is successfully completed only when $k$ identical examinations are successfully completed.

Sandhya *et al*. (2005) have shown the Harris distribution as a gamma mixture by considering a linear function of Poisson r.v with gamma distributed Poisson parameter. Also a characterization of exponential distribution using a Laplace transform was given in Sandhya and Satheesh (1997). Motivated by this characterization gamma distribution is characterized in section 4.

## 2. Extended geometric distribution and its mixture

We prove the following characterizations of extended geometric distribution.

*Theorem 1:* A probability sequence (PS) $\{f(a+nk)\}$, $n = 0,1,2, \ldots$; $a \geq 0$ and $k > 0$ are integers is completely monotone if and only if it is a mixture of extended geometric distribution.

*Proof:* Consider an extended geometric r.v with PS given by (1). Then we have

$$P(Y < a+nk) = 1 - P(Y \geq a+nk) = 1 - q^n$$

Randomizing $q$ with a distribution $G$ concentrated on (0,1), we have the distribution function (DF) $F$ of the mixture as

$$F(a+nk) = P(X < a+nk) = 1 - \int_0^1 q^n dG(q)$$
$$F(a+nk) = 1 - m(n), \quad n = 0,1,2,\ldots$$

where $\{m(n)\}$, the moment sequence of $G$ is completely monotone (CM) by Hausdorff's theorem given in Feller (1966, p. 223). If $\{f(a+nk)\}$ denote the PS corresponding to this mixture, then,

$$\Delta F(a+nk) = f(a+nk) = -\Delta m(n), \quad n = 0,1,2,\ldots$$

and $\{f(a+nk)\}$ is CM, $\Delta$ being the differencing operator.

Conversely, starting with a CMPS $\{f(a+nk)\}$, the survival sequence $S(a+nk) = P(X \geq a+nk)$ is again CM. Also, when $n = 0$, $S(a+nk) = S(a) = P(X \geq a) = 1$. Now by Hausdorff's theorem and retracing the steps we see that $\{f(a+nk)\}$ corresponds to a mixture of extended geometric distribution and the proof is complete.

Now consider a non-negative integer-valued r.v $U$. Then the r.v $V = kU$, $k > 0$ an integer, is with PS $\{g(nk)\}$, where $g(nk) = P(V = nk)$, $n = 0,1,2,\ldots$. The survival sequence is taken as

$$S(nk) = g((n+1)k) + g((n+2)k) + \ldots$$

*Theorem 2:* A PS $\{g(nk)\}$, $k > 0$ an integer satisfies the convolution equation

$$\{g(nk)\} * \{S(nk)\} = \{ng(nk)\}, \quad n = 0,1,2,\ldots \text{ if and only if}$$

$$g(nk) = q^{n-1}p, \quad n = 1,2,\ldots, \quad 0 < p < 1, \quad p + q = 1.$$

*Proof:* The survival sequence corresponding to the PS $g(nk) = q^{n-1}p$, $n = 1, 2, \ldots$ is given by $S(nk) = q^n$, $n = 0, 1, 2, \ldots$, which satisfies the convolution equation

$$\{g(nk)\} * \{S(nk)\} = \{ng(nk)\}, \quad n = 0,1,2, \ldots$$

Now suppose that the convolution equation is true. Then taking the corresponding p.g.f's and using the theorem given in Feller (1968, p.267) we have

$$P(s) \cdot \frac{1 - P(s)}{1 - s^k} = \frac{s}{k} \cdot P'(s)$$

where $P(s)$ is the p.g.f of $\{g(nk)\}$. The solution of this differential equation for a p.g.f is

$$P(s) = \frac{as^k}{a + b - bs^k} = \frac{ps^k}{1 - qs^k}, \quad p = \frac{a}{a+b}, \quad q = 1 - p$$

This is the p.g.f of a $Geo_k(p, k)$ distribution.

*Theorem 3:* A PS $\{g(nk)\}$, $k > 0$ an integer satisfies the convolution equation

$$\{g(nk)\} * \{S(nk)\} = \{(n+1)g((n+1)k)\}, \quad n = 0,1,2,\ldots \text{ if and only if}$$

$$g(nk) = q^n p, \quad n = 0,1,2,\ldots, \quad 0 < p < 1, \quad p + q = 1.$$

*Proof:* Corresponding to the PS $g(nk) = q^n p$, $n = 0,1,2,\ldots$ we have the survival sequence $S(nk) = q^{n+1}$, $n = 0,1,2,\ldots$. Then it can be shown that

$$\{g(nk)\} * \{S(nk)\} = \{(n+1)g((n+1)k)\}, \quad n = 0,1,2,\ldots$$

Now, assume that the convolution equation is true. Then proceeding as in *theorem 2* we have:

$$P(s) \cdot \frac{1-P(s)}{1-s^k} = \frac{P'(s)}{k s^{k-1}}$$

where $P(s)$ is the p.g.f of $\{g(nk)\}$. The solution of this differential equation for a p.g.f is

$$P(s) = \frac{c}{1+c-s^k} = \frac{p}{1-qs^k}, \quad p = \frac{c}{1+c}, \quad q = 1-p.$$

This is the p.g.f of a $Geo_0$ $(p, k)$ distribution.

### 3. Characterizations of Harris and negative binomial distributions

In this section we characterize the Harris distribution $H_1(m,k,1/k)$. Its p.g.f $P(s)$ given in (2) on differentiation with respect to $s$ gives,

$$P'(s) = m/\left(m-(m-1)s^k\right)^{1+(1/k)}, \text{ where } m = P'(1) \tag{4}$$

Now we have the following theorems.

*Theorem 4:* A r.v $X \sim H_1(m,k,1/k)$ if and only if its p.g.f $P(s)$ satisfies the equation

$$s^{1+k} P'(s) = m(P(s))^{1+k}$$

*Proof:* Let $X \sim H_1(m,k,1/k)$. Then we have (2) and (4) which on simplification gives

$$s^{1+k} P'(s) = m(P(s))^{1+k} \tag{5}$$

Conversely, starting with (5) we have,

$$\int dP/P^{1+k} = m \int ds/s^{1+k} \quad \text{and hence} \quad 1/P^k = (m/s^k) - c$$

When $s = 1$, $P(s) = 1$ and hence $c = m-1$. So $P(s) = s/\left(m-(m-1)s^k\right)^{1/k}$, which is the p.g.f of the Harris distribution $H_1(m,k,1/k)$.

*Theorem 5:* A Harris distribution $H_1(m,k,1/k)$ is characterized by the equation

$$s(1-s^k) P'(s) = P(s)\{1- (P(s))^k\}.$$

*Proof:* Let $X \sim H_1(m,k,1/k)$. Using (2) and (4) we get (5). Again, from (2) we have

$$m - (m-1) s^k = (s/P(s))^k \quad \text{which gives}$$

$$m = (s/P(s))^k \left\{\left(1-(P(s))^k\right)/\left(1-s^k\right)\right\}.$$

Now substituting this $m$ in (5) we get

$$s(1-s^k)P'(s) = P(s)\{1-(P(s))^k\} \tag{6}$$

Conversely, starting with a p.g.f satisfying (6) we have

$$\int dP / P(1-P^k) = \int ds / s(1-s^k).$$

Now, this on integration and simplification gives $P^k/(1-P^k) = cs^k/(1-s^k)$. Then the solution $P(s) = s/\left(m-(m-1)s^k\right)^{1/k}$, $m=1/c$, is the p.g.f of the Harris distribution $H_1(m,k,1/k)$.

*Corollary:* For a Harris distribution $H_1(m,k,1/k)$, $k$ may be evaluated as

$$k = \frac{sP(s)P''(s) - s(P'(s))^2 + P(s)P'(s)}{P'(s)(sP'(s) - P(s))}$$

where $P'(s)$ and $P''(s)$ are the first and second order derivatives of $P(s)$ with respect to $s$.

Proceeding as given in *theorem 5* we have,

$$s(1-s^k)P'(s) = P(s)\{1-(P(s))^k\}.$$

Now from (2) and (4) we get, $P'(s) = m(P(s)/s)^{1+k}$ which on differentiation with respect to $s$ gives $P''(s) = m(1+k)(P(s))^k\{sP'(s) - P(s)\}/s^{k+2}$. Using this we get

$$\frac{P''(s)}{P'(s)} = \frac{(1+k)\{sP'(s) - P(s)\}}{sP(s)}, \text{ which gives } k = \frac{sP(s)P''(s) - s(P'(s))^2 + P(s)P'(s)}{P'(s)(sP'(s) - P(s))}.$$

Next, we have a theorem which shows the relationship between a Harris distribution and an extended geometric distribution.

*Theorem 6:* $X_1, X_2, \ldots, X_k$ are $k$ i.i.d $H_1(m,k,1/k)$ r.v s if and only if $Y = \sum_{i=1}^{k} X_i$ is an extended geometric random variable.

*Proof:* Let $X_1, X_2, \ldots, X_k$ be $k$ i.i.d $H_1(m,k,1/k)$ r.v s with p.g.f

$$P_{X_i}(s) = \frac{s}{\left(m-(m-1)s^k\right)^{1/k}}, \quad i = 1, 2, \ldots, k.$$

The p.g.f of $Y = \sum_{i=1}^{k} X_i$ is then

$$P_Y(s) = (P_{X_i}(s))^k = \frac{s^k}{m-(m-1)s^k}$$

$$= \frac{ps^k}{1-qs^k}, \quad p = 1/m, \, q = 1-p, \, 0<p<1$$

This is the p.g.f of a $Geo_k(p, k)$ distribution. Hence $Y = \sum_{i=1}^{k} X_i \sim Geo_k(p,k)$.

Conversely, let $Y = \sum_{i=1}^{k} X_i \sim Geo_k(p,k)$ where $X_i$s are i.i.d. Now, the p.g.f of $Y$ is given by

$$P_Y(s) = \frac{ps^k}{1-qs^k} = \frac{s^k}{m-(m-1)s^k} = \left(\frac{s}{(m-(m-1)s^k)^{1/k}}\right)^k = (P_{X_i}(s))^k, \quad m = 1/p$$

We can see that $P_{X_i}(s)$ is the p.g.f of a $H_1(m,k,1/k)$ r.v and hence $X_1, X_2, \ldots, X_k$ are $k$ i.i.d $H_1(m,k,1/k)$ r.v s.

Arguing on similar lines we have:

*Theorem 7:* $X_1, X_2, \ldots, X_k$ are $k$ i.i.d negative binomial $NB(p, 1/k)$ variables if and only if $Y = \sum_{i=1}^{k} X_i$ is a geometric random variable.

For some connected results see Satheesh *et al.* (2005, available at *http://arxiv.org/abs//math.PR/0507535*).

Again, the p.g.f of a $H_1(m,k,1/k)$ distribution can be written as

$$P(s) = \frac{s}{(m-(m-1)s^k)^{1/k}} = \left(\frac{s^k}{m-(m-1)s^k}\right)^{1/k} = (P_1(s))^{1/k}$$

where $P_1(s) = \frac{s^k}{m-(m-1)s^k} = \frac{ps^k}{1-qs^k}$, $p = 1/m$ is the p.g.f of $Geo_k(p, k)$ distribution.

Thus the $k^{th}$ root of the p.g.f of an extended geometric distribution defined on $\{k, 2k, \ldots\}$ is the p.g.f of the Harris distribution.

## 4. Characterization of gamma distribution

Here we present a characterization of gamma distribution.

*Theorem 8:* The Laplace transform $L(s)$ of a probability distribution satisfies the differential equation $L'(s) = L(s)\log L(s) /(1+s)\log(1+s)$ if and only if the distribution is gamma$(\alpha,1)$ with density function $f(x) = e^{-x} x^{\alpha-1}/\Gamma\alpha$, $x > 0$, $\alpha > 0$.

*Proof:* Let $X$ be a gamma$(\alpha,1)$ r.v with Laplace transform $L(s) = (1+s)^{-\alpha}$. Hence $\alpha = -\log L(s) / \log (1+s)$. Differentiation of $L(s)$ with respect to $s$ gives

$$L'(s) = -\alpha L(s)/(1+s)$$

Substituting the value of $\alpha$ we get

$$L'(s) = L(s) \log L(s)/ (1+s) \log(1+s) \tag{7}$$

Now, suppose that equation (7) is true. Then we have,

$$\int dL(s)/L(s)\log L(s) = \int ds/(1+s)\log(1+s)$$

$u = \log L(s)$ and $t = \log(1+s)$ gives $u = ct$. i.e. $L(s) = (1+s)^c$ which gives $c = \log L(s)/\log(1+s)$ and hence from (7), $L'(s) = c L(s)/(1+s)$. When $s = 0$, $L(0) = 1$ and $L'(0) = -\alpha$ (say) and hence $c = -\alpha$. Hence $L(s) = (1+s)^{-\alpha}$, which is the Laplace transform of a gamma distribution with parameter $\alpha$.